\documentclass[journal]{IEEEtran}
\usepackage{cite}
\usepackage{graphicx}
\usepackage{url}
\usepackage[colorlinks=true,linkcolor=black,anchorcolor=black,citecolor=black,urlcolor=black]{hyperref}
\usepackage{lipsum}
\usepackage{color}
\usepackage{amsmath}
\usepackage{algorithmic}
\usepackage{array}
\usepackage{multirow}
\usepackage{makecell}
\usepackage{float}
\usepackage{setspace}
\usepackage{array}
\usepackage{supertabular}
\usepackage{supertabular,booktabs}
\usepackage{amssymb}
\usepackage{nomencl}
\makenomenclature

\usepackage{amsfonts}
\usepackage{comment}
\newboolean{showcomments}
\setboolean{showcomments}{true}
\newcommand{\wwei}[1]{\ifthenelse{\boolean{showcomments}}
	{ \textcolor{red}{(WW:  #1)}}{}}
\newcommand{\bliu}[1]{\ifthenelse{\boolean{showcomments}}
	{ \textcolor{blue}{(BL:  #1)}}{}}

\ifCLASSOPTIONcompsoc
\usepackage[caption=false,font=normalsize,labelfont=sf,textfont=sf]{subfig}
\else
\usepackage[caption=false,font=footnotesize]{subfig}
\fi
\usepackage{fixltx2e}
\usepackage{url}
\usepackage{amsfonts}
\usepackage{multicol}
\begin{document}
	\title{\huge Unbalance Mitigation via Phase-switching Device and Static Var Compensator in Low-voltage Distribution Network}

%	\author{\thanks{Manuscript received XX; revised XX. ({\it Corresponding author: Bin Liu}). This work is funded in part by ARENA Project: Demonstration of three-phase dynamic grid-side technologies for increasing distribution network DER hosting capacity.}}

	\author{Bin Liu, {\it Member, IEEE}, Ke Meng, {\it Senior Member, IEEE}, Zhao Yang Dong, {\it Fellow, IEEE}, Peter K. C. Wong, {\it Senior Member, IEEE}, Tian Ting
		\thanks{Manuscript received XX; revised XX. ({\it Corresponding author: Bin Liu}). This work is funded by ARENA (Australian Renewable Energy Agency) Project: Demonstration of three-phase dynamic grid-side technologies for increasing distribution network DER hosting capacity.}
		\thanks{Bin Liu, Ke Meng and Zhao Yang Dong are with School of Electrical Engineering and Telecommunications, The University of New South Wales, Sydney 2052, Australia. (e-mail: bin.liu@unsw.edu.au, ke.meng@unsw.edu.au, zydong@ieee.org)}
		\thanks{Peter K. C. Wong is with Jemena Electricity Networks (Vic) Ltd, Melbourne 3000, Australia. (e-mail: peter.wong@jemena.com.au)}
		\thanks{Tian Ting is with AusNet Services, Melbourne 3000, Australia. (e-mail: tian.ting@ausnetservices.com.au)}
	}

	\markboth{\LaTeX~Journal Paper,~Vol.~XX, No.~XX, Aug~2019}%
	{Shell \MakeLowercase{\textit{et al.}}: Bare Demo of IEEEtran.cls for IEEE Journals}

	\maketitle

	\begin{abstract}
As rooftop solar PVs installed by residential customers penetrate in low voltage distribution network (LVDN), some issues, e.g. over/under voltage and unbalances, which may undermine the network's operational performance, need to be effectively addressed. To mitigate unbalances in LVDN, dynamic switching devices (PSDs) and static var compensator (SVC) are two equipment that are cost-effective and efficient. However, most existing research on operating PSDs are based on inflexible heuristic algorithms or without considering the network formulation, which may lead to strategies that violate operational requirements. Moreover, few pieces of literature have been reported on mitigating unbalances in LVDN via SVC and PSDs together. This paper, after presenting the dispatch model of SVC, formulates the decision-making process as a mixed-integer non-convex programming (MINCP) problem considering all practical operational requirements. To efficiently solve the challenging problem, the MINCP is reformulated as a mixed-integer second order-cone programming (MISOCP) problem based on either exact reformulations or accurate approximations, making it possible to employ efficient off-the-shelf solvers. Simulations based on a modified IEEE system and a practical system in Australia demonstrates the efficiency of the proposed method in mitigating unbalances in LVDN.     
	\end{abstract}
	% \fi
	\IEEEpeerreviewmaketitle
	
	%Section I Introduction
	\section{Introduction}
    \IEEEPARstart{T}{o} address the sustainability and environmental issue, renewable energy has been developing rapidly worldwide and is expected to be continuously growing in the future. Due to encouraging policies and incentives from government agencies, Australia is experiencing a remarkable renewable energy development in recent years, taking the top spot worldwide in the penetration level of residential PV installation in the low-voltage distribution network (LVDN). The annual report released by Australian Energy Council shows that the year 2018 was a record-breaking year for PV development in Australia, with total installed residential capacity reaching over 1.4 GW, which increases by 20\% compared with the year 2017. By the end of the year 2018, the cumulative installed capacity of residential PV in Australia stood at 7.98 GW with more than 2 million installations across the nation, and the numbers keep growing \cite{RN192}. 
	
    Like European LVDN, which is different from that in North America, power utilities in Australia run extensive four-wire (230/400V) grid along the streetscape \cite{RN61}, where most residential customers are powered by a single phase cable. As the number of installed PVs continues increasing, unbalance issue may be worsened during the hours with high solar radiation level\cite{RN52}. This is because most LVDNs are not specially designed to accommodate the distributed generation sources in the design stage. Even if PVs are carefully placed in the network, the diversity of electricity usage profiles of all customers can still cause significant unbalances in the network. Unbalances, either current unbalance or voltage unbalance, in LVDN may lead to lots of operational problems, e.g. increased power loss caused by high neutral current, contributing to over-/under-voltage issue in one phase due to light/stressed power flow and shortened lives of appliances \cite{RN103,RN104}. 
	
\begin{figure}[htb!]
		\centering\includegraphics[scale=0.55]{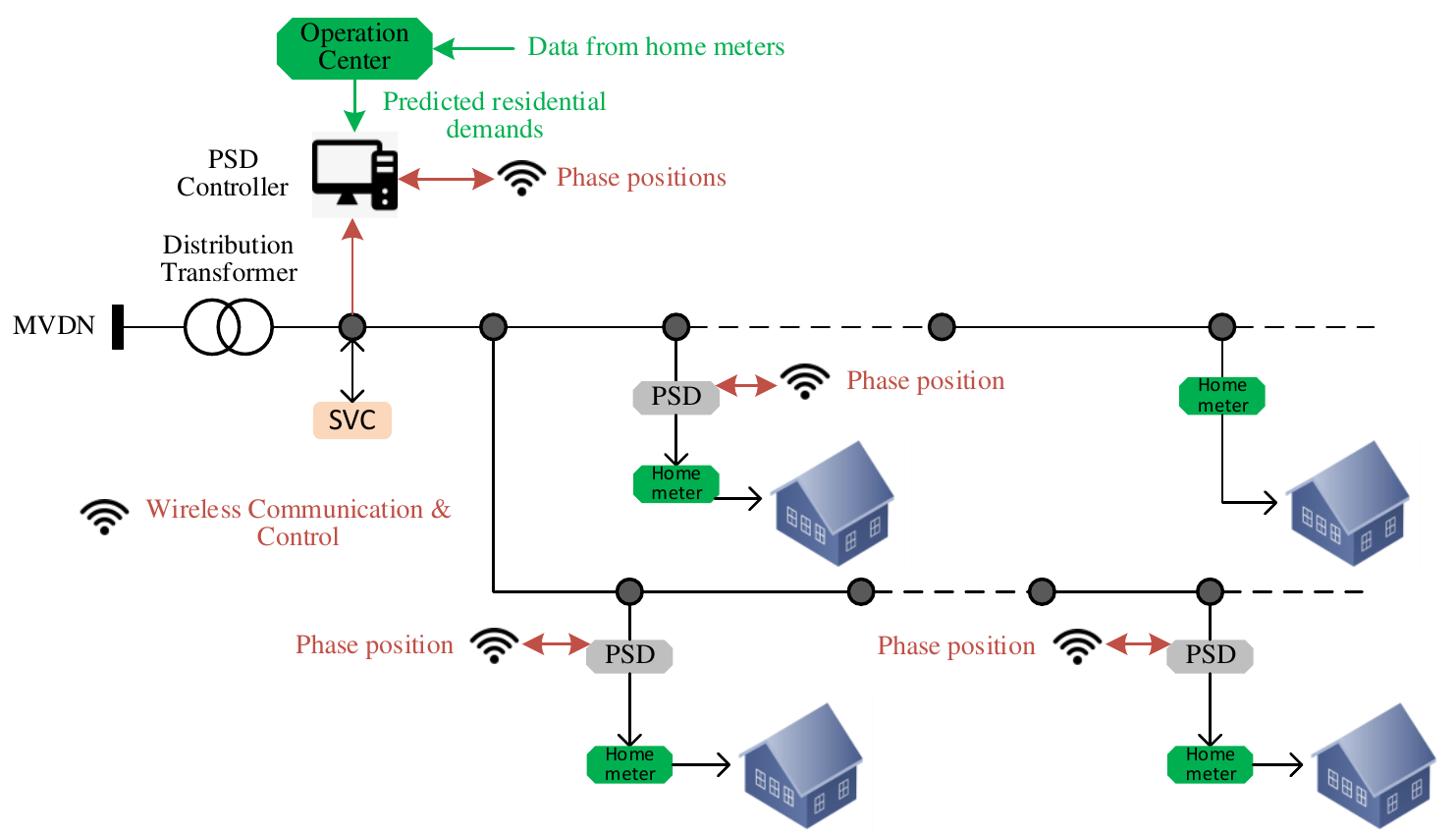}
		\caption{Communication and control system for controlling PSDs and SVC in LVDN (MVDN: Medium voltage distribution network)}
		\label{fig-PSDController-Smeter-Operational}
	\end{figure}

    To mitigate the unbalances, a cost-efficient option is using the phase-switching device (PSD) to change the connected phase of each residential customer regularly or when monitored unbalance exceeds a predefined level \cite{RN61,RN44}. Benefiting from developments of communication and control technologies, existing mechanical PSD, which can merely be adjusted by cost-intensive labors, can be replaced by smart PSD that is manufactured with modern electronic technologies to achieve a smooth switching process and minimal impacts on residential customers. As shown in Fig.\ref{fig-PSDController-Smeter-Operational}, with predicted residential customer demands, the controller, which is usually mounted at the secondary side of the distribution transformer (DT), can optimize the phase positions of all PSDs and send the control signals to each device through the wireless network.
	%Operating PSD in LVDN is actually a three-phase unbalanced optimal power flow (TUOPF) problem with introduced integer variables indicating the positions of PSDs. 
	%Reported methods to solve TUOPF can be divided into three categories: iteration-based method \cite{RN35,RN18,RN42,RN63}, LP-based method \cite{RN50,RN66,RN62,RN67,RN68,RN41,RN30} and convex relaxation-based method \cite{RN54,RN50,RN53,RN45,RN40,RN64}.
		  
    Reported methods to address unbalance issue in LVDN via PSDs are mostly based on heuristic algorithm supported by independent three-phase unbalance power flow (TUPF) programs \cite{RN104,RN99,RN71,RN77,RN98,RN74,RN76,RN101,RN106}. A feasible solution can be provided by this type of method after a reasonable time. However, the computation efficiency may become an issue with a large number of introduced integer variables. Moreover, the operational constraints, e.g. voltage magnitude limits, current magnitude limits and voltage unbalance level, are rarely considered. To make the formulation more flexible and improve the computational efficiency, linear or convex models are proposed based on various simplifications. For example, a mixed-integer linear programming (MILP) formulation was presented in \cite{RN103} by treating all customers as constant-current loads. A linear programming (LP) formulation is presented in \cite{RN44}, and the robust method is presented in \cite{RN105} to further take load profile uncertainty into account. However, neither \cite{RN44} nor \cite{RN105} considers the network formulation when seeking the optimal strategy. 
	
    %However, PSDs are only capable of providing rough regulation in the network, as we explain through the following illustrative example. Assuming there are 4 resistive customers in an LVDN and they are connected to phase $a,b,c,a$ with net demands being 8 kW, 5 kW, 2 kW  and 3 kW, respectively. If only customer 4 is quipped with PSD, then the best strategy is switching it from phase $a$ to phase $c$ to minimize power difference between any two phases. Obviously, 3 kW of power unbalance still exists, implying that the system can only be roughly balanced by PSD. 
	
    With the development of power electronic technologies, unbalances in the network can be further regulated by operating PSD and static var compensator (SVC) together to achieve better performance. Traditionally, the network can be balanced by passive equipment such as capacitors and inductors based on Steinmetz circuit \cite{SVGbook,SVGpaper02,4039441,SVC-1,SVC-2,SVC-3}. Each susceptance connected between unbalanced phases in that case is usually fixed, which lacks flexibility because loads of three phases may vary significantly at different periods. However, with the development of power electronic technologies, a variable susceptance can be achieved by using controlled thyristor, i.e. via SVC.

	%which can not only compensate reactive power to the network, but also transferring active power among three phases, thus is capable of achieving better operation performance of LVDN in terms of power balances \cite{SVGbook,SVGpaper01,SVGpaper02,4039441}. For example, 2 kW of demand can be distributed from phase $a$ to phase $b$ and $c$ equally at the secondary side of DT by SVG in the illustrative example to realize perfect power balance running through the DT. In other words, the LVDN is balanced if looked at from the upstream network. 
	
	%With the development of power electronics technology, unbalances in the network can be further regulated with higher granularity. One of such equipment is Static VAR Generator (SVG), which can not only compensate reactive power to the network, but also transferring active power among three phases, thus is capable of achieving better operation performance of LVDN in terms of power balances \cite{SVGbook,SVGpaper01,SVGpaper02,4039441}. For example, 2 kW of demand can be distributed from phase $a$ to phase $b$ and $c$ equally at the secondary side of DT by SVG in the illustrative example to realize perfect power balance running through the DT. In other words, the LVDN is balanced if looked at from the upstream network. 
		
    To the authors' knowledge, few pieces of literature have been reported on mitigating unbalances via PSD and SVC together in LVDN, which motivates us to investigate flexible and efficient methods in this area while taking all operational constraints into account. %It is expected that the voltage levels and power flows in the network can be ameliorated via coordinately operating the two devices. 
	 
	%It is noteworthy that by balancing three phases via SVG, which can affect the voltage level of the whole network, and optimally setting phase positions of all PSDs, which can help ameliorate power flows and voltage levels in the network, the unbalance issue can be well addressed while ensuring the network operates within practical requirements.
	
	%\textcolor{red}{find two cases: (1) SVG is not feasible while will be feasible with SVG\&PSD; (2) PSD is not feasible while with SVG\&PSD is feasible}
	
    Compared with existing research, the contributions of this paper are summarized as follows. 
	\begin{enumerate}
		\item An efficient SVC model for dispatch purposes in LVDN is presented. The formulation takes the controllability of SVC into account and can realize its functionalities of both compensating reactive currents and transferring active currents among three phases. 
		\item The mixed-integer non-convex programming (MINCP) formulation of operating PSD and SVC together in LVDN to mitigate unbalances is presented based on an efficient linear power flow model. The formulation is flexible and takes practical operational constraints into account.
		\item The non-convexity in the formulated MINCP problem is systematically studied, which is thereafter reformulated to a mixed-integer second-order cone (MISOCP) problem based on several reasonable and practical assumptions, thus making it possible to employ commercial solvers to efficiently solve the challenging problem.
	\end{enumerate}
	
     The remainder of this paper is organized as follows. The model of SVC for dispatch purposes in LVDN is presented in Section II, followed by the optimization model to coordinately operate PSDs and SVC in LVDN. The non-convexity analysis and solution techniques are presented in Section III. Case studies based on a modified IEEE system and a practical system in Australia are performed in Section IV and the paper is concluded in Section V.
	
%=====================================================================================================================================
\section{Problem Formulation}
In the formulation, $\mathcal {C}_i$ will be used to represent the set of customers connected to node $i$, $\mathcal{F}_i$ and $\mathcal{X}_i$ as the sets of customers without and with PSDs installed (denoted as fixed and adjustable customers), respectively. Therefore, there are $\mathcal{F}_i\cap \mathcal{X}_i=\emptyset$ and $\mathcal{F}_i\cup \mathcal{X}_i=\mathcal{C}_i$. Moreover, $V$ and $U$ represent the voltage of a node in the main feeder and terminal voltages of a customer, respectively. $X$ and $Y$ are the real and imaginary parts of $V$; $I$ represents the current with $J$ and $W$ being the real and imaginary parts, respectively; The subscript $\phi/\psi, t,i,j,ik$ represent phase $\phi/\psi$ that belongs to $\{a,b,c\}$, period $t$, node $i$, the $j^\text{th}$ customer at node $i$ and line $ik$, respectively. Other parameters or variables will be explained right after their appearances.

%Several assumptions are listed below based on practical operational information in Australia.
%\begin{enumerate}
%	\item The voltage of root node, i.e. the primary side of DT, is known noting that this voltage is mainly determined by the medium-voltage distribution network (MVDN) in upper stream.
%	\item The residential loads for all customers are with constant PQ values in each period.
%	\item Feeders between any two poles in LVDN are constructed with four-wire (phase $a,b,c$ and zero earthed conductor), which is the general case in Australia \cite{RN61}.
%\end{enumerate}

\subsection{Dispatch model of SVC in LVDN}
\begin{figure}[htb!]
	\centering\includegraphics[scale=0.46]{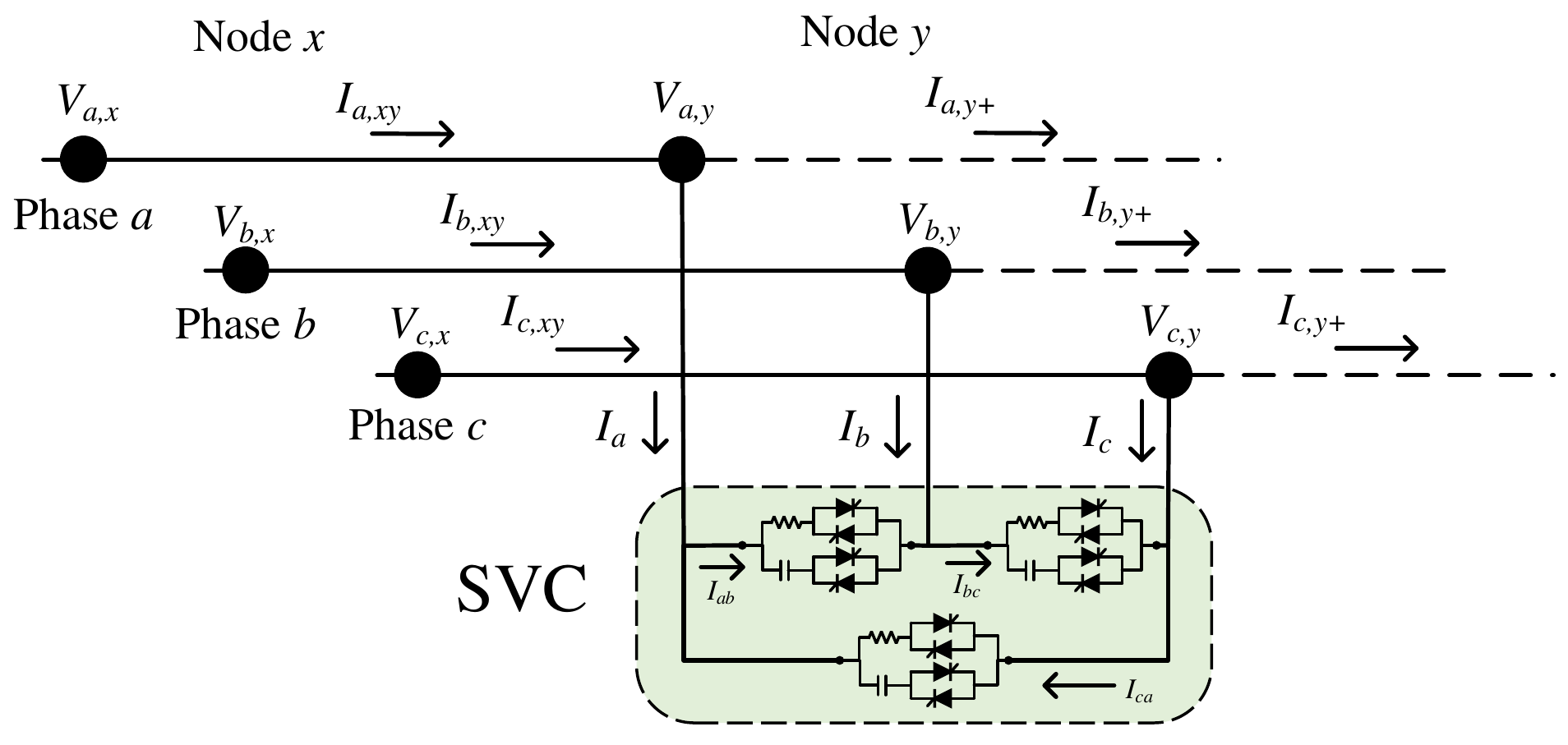}
    \caption{Connection topology of SVC in LVDN ($I_{ab},I_{bc}$ and $I_{ca}$ are control variables)}
	\label{fig-svgModel-01}
\end{figure}

Assuming the node indices of the primary and secondary sides of the DT are $x$ and $y$, respectively, the dispatch model of SVC can be expressed as follows based on Fig.\ref{fig-svgModel-01}.
\begin{subequations}
	\small
	\label{Ipcd}
	\begin{eqnarray}
	\label{Ipcd-1}
	I_{a,t}=I_{ab,t}-I_{ca,t}~\forall t\\
	\label{Ipcd-2}
	I_{b,t}=I_{bc,t}-I_{ab,t}~\forall t\\
	\label{Ipcd-3}
	I_{c,t}=I_{ca,t}-I_{bc,t}~\forall t\\
	\label{Ipcd-4}
	I_{\phi\psi,t}=|I_{\phi\psi,t}|(\cos{\beta_{\phi\psi,y,t}}+j\sin{\beta_{\phi\psi,y,t}})~\forall\phi\psi\in\Phi,\forall t\\
	\label{Ipcd-5}
	|I_{\phi\psi,t}|\le \frac{\kappa_{\phi\psi,t} S^{cap}+(1-\kappa_{\phi\psi,t}) S^{ind}}{3|V^{svc}_{n}|}~\forall\phi\psi\in\Phi,\forall t\\
	\label{Ipcd-6}
	\beta_{\phi\psi,t}=\delta _{\phi\psi,y,t}+\kappa_{\phi\psi,t}\pi-\pi/2~\forall\phi\psi\in\Phi,\forall t\\
	\label{Ipcd-7}
	\kappa_{\phi\psi,t}\in\{0,1\}~\forall t    
	\end{eqnarray}
\end{subequations}
where $\delta_{\phi,y,t}$ and $\delta_{\phi\psi,y,t}$ are the phase angles of $V_{\phi,y,t}$ and $V_{\phi\psi,y,t}$, respectively; 
$I_{\phi,t}$ is the load current of SVC at period $t$; 
$I_{\phi\psi,t}$ is the current running in SVC flowing from phase $\phi$ to phase $\psi$ at period $t$;
$\kappa_{\phi\psi,t}$ is an introduced binary variable indicating whether the impedance between phase $\phi$ and phase $\psi$ is inductive ($\kappa_{\phi\psi,t}=0$) or capacitive ($\kappa_{\phi\psi,t}=1$) at period $t$;
$V^{svc}_{n}$ is the rated phase-to-phase voltage of the SVC;
$S^{cap}$ and $S^{ind}$ are the capacitive and inductive capacity of the SVC; 
$\beta_{\phi\psi,t}$ is the phase angle of $I_{\phi\psi,t}$ and 
$\Phi=\{ab,bc,ca\}$

The load currents of SVC in three phases are calculated based on \eqref{Ipcd-1}-\eqref{Ipcd-4} according to Kirchhoff's current law (KCL). The magnitudes of phase-to-phase currents provided by SVC are constrained by \eqref{Ipcd-5}. Obviously, the dispatch model of SVC is strongly non-convex due to \eqref{Ipcd-4}, the approximation of which will be further discussed in the next section. \eqref{Ipcd-6}-\eqref{Ipcd-7} are to make sure the phase angle of $I_{\phi\psi,t}$ is always 90$^\circ$ lagging or leading $V_{\phi\psi,y,t}$.
\begin{figure}[htb!]
	\centering\includegraphics[scale=0.17]{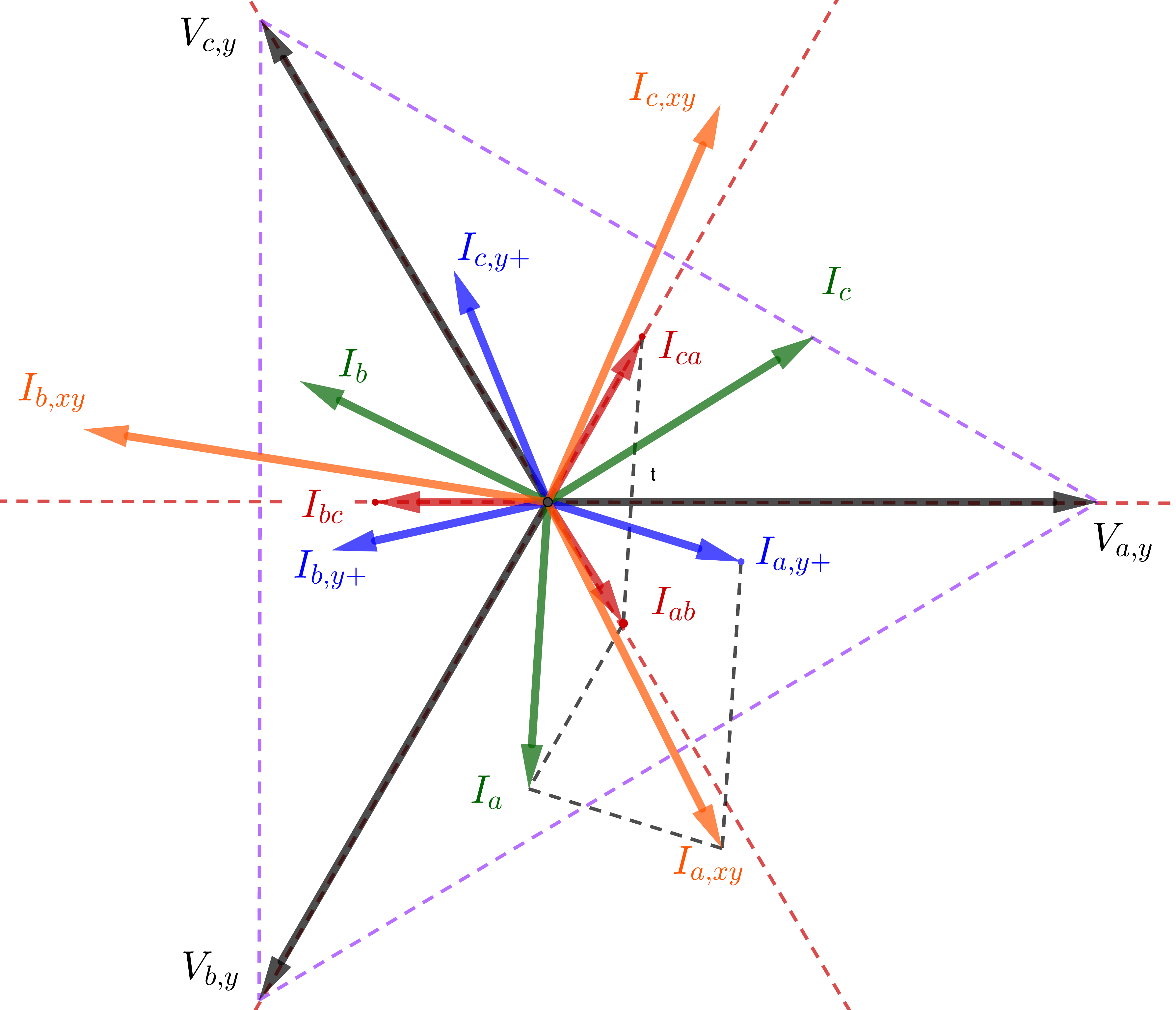}
	\caption{Phasor illustration of current compensation from SVC ($I_{\phi,y+}$ is the total current in the downstream of node $y$ when SVC is disconnected)}
	\label{fig-svgModel-02}
\end{figure}

The basic idea of balancing currents via SVC can be briefly illustrated by Fig.\ref{fig-svgModel-02}, where $I_{ab},I_{bc}$ and $I_{ca}$ are control variables and are always perpendicular to $V_{a,y}-V_{b,y},V_{b,y}-V_{c,y}$ and $V_{c,y}-V_{a,y}$, respectively. In the figure, $I_{a,y+},I_{b,y+}$ and $I_{c,y+}$, which are the currents running through the DT when SVC is disconnected, are initially with strong unbalance. However, by controlling $I_{\phi\psi}(\forall \phi\psi\in\{ab,bc,ca\})$, $I_{a},I_{b}$ and $I_c$ can be regulated and, if with appropriate control, lead to much more balanced currents running through the DT, which equals to $I_{\phi,xy}=I_{\phi}+I_{\phi,y+}(\forall \phi)$. Noting that there may be $V_{\phi,y}I^*_{\phi,y+}\neq V_{\phi,y}I^*_{\phi,xy}$, which implies that through controlling SVC reactive powers in three phases running through the DT can be regulated while active powers can be redistributed among three phases.

\subsection{Objective}
In this paper, the objective is to minimize the unbalance level of currents in three phases running through the DT. In this paper, the  negative sequence and zero-sequence components of the currents in three phases are considered, and the objective is accordingly formulated as the sum of magnitudes of them in all concerned dispatching periods, yielding 
\begin{eqnarray}
\label{obj-1-opsd}
F=\sum\nolimits_{t\in\cup_k\mathcal{T}_k}{z_t}=\sum\nolimits_{t\in\cup_k\mathcal{T}_k}{(z^-_{t}+z^0_{t})}
\end{eqnarray}
where $\mathcal{T}_k$ is the $k^\text{th}$ set of concerned dispatch periods;
% $xy$ is the index of root line, which is the line starting from the secondary side of DT; 
$z^-_{t}$ and $z^0_{t}$ are the magnitudes of negative sequence and zero-sequence currents running through the DT at period $t$, respectively.

Denoting the current of phase $\phi$ running through the DT at period $t$ as $I_{\phi,xy,t}$ with subscript $xy$ indicating the line where the DT locates, $z^-_{t}$ and $z^0_{t}$ can be formulated as follows.
\begin{subequations}
	\small
	\label{ulevel-PCD-opsd}
	\begin{eqnarray}
	\label{ulevel-PCD-opsd-1}
	|I^-_{t}|\le z_t^-,~|I^0_{t}|\le z_t^0~\forall t\\
	\label{ulevel-PCD-opsd-2}
	6I^-_{t}=2I_{a,xy,t}-(1+j\sqrt{3})I_{b,xy,t}-(1-j\sqrt{3})I_{c,xy,t}~\forall t\\
	\label{ulevel-PCD-opsd-3}
	3I^0_{t}=I_{a,xy,t}+I_{b,xy,t}+I_{c,xy,t}~\forall t
	\end{eqnarray}
\end{subequations}
where $I^-_{t}$ and $I^0_{t}$ are the negative sequence and zero-sequence currents running through the DT in phase $a$ at period $t$, respectively.

It is noteworthy that as both the magnitudes of negative sequence and zero-sequence currents running through the DT are identical for three phases,  \eqref{ulevel-PCD-opsd-2} and \eqref{ulevel-PCD-opsd-3} are constructed based on phase $a$ only.

\subsection{Operational constraints}
   \subsubsection{Ohm's law for each line in main feeder}
    The Ohm's law for each line in main feeder is formulated as
	\begin{eqnarray}
	\label{ol-1-opsd}
	V_{\phi,i,t}-V_{\phi,k,t}=\sum\nolimits_\psi{Z_{i,k}^{\phi\psi}I_{\phi,ik,t}}~\forall \phi,\forall ik,\forall t
	\end{eqnarray}
	where 
	%$V_{\phi,i,t}$ and $I_{\phi,ik,t}$ are the voltage of node $i$ and current of line $ik$ for period $t$ in phase $\phi$, respectively;
	$Z_{i,k}^{\phi\psi}$ is the mutual impedance between phase $\phi$ and $\psi$ of line $ik$.
	
	\subsubsection{Kirchhoff's current law (KCL) at each node}
	This constraint ensures current balance at each node as shown in \eqref{kcl-1-opsd} and \eqref{kcl-1-opsd+}. Customer's demand current and the current flowing out of its connected node are bridged by \eqref{kcl-2-opsd}. Moreover, the customer can only be connected to one phase as constrained by \eqref{kcl-4-opsd} and \eqref{kcl-5-opsd}.
	\begin{subequations}
		\small
		\label{kcl-opsd1}
		\begin{eqnarray}
		\label{kcl-1-opsd}
		\sum_{n:n\rightarrow i}{I_{\phi,ni,t}}-\sum_{k:i\rightarrow k}{I_{\phi,ik,t}}=\sum_{j\in \mathcal {C}_i}{I_{\phi,i,j,t}}~\forall \phi,\forall i\notin\{x,y\},\forall t\\
		\label{kcl-1-opsd+}
		I_{\phi,xy,t}-\sum_{k:y\rightarrow k}{I_{\phi,yk,t}}-I_{\phi,t}=\sum_{j\in \mathcal {C}_y}{I_{\phi,y,j,t}}~\forall \phi,\forall t\\
		\label{kcl-2-opsd}
		I_{\phi,i,j,t}=\varepsilon_{\phi,i,j,t}I_{i,j,t}=
		\left\{\begin{array}{c}
			\alpha_{\phi,i,j,t}I_{i,j,t}~\forall j\in\mathcal{F}_i\\
			\mu_{\phi,i,j}I_{i,j,t}~\forall j\in\mathcal{X}_i
		\end{array}\right.	~\forall \phi,\forall i,\forall t\\
%		I_{\phi,i,j,t}=\alpha_{\phi,i,j,t}I_{i,j,t}~\forall \phi,\forall i,\forall j\in \mathcal{F}_i,\forall t\\
%		\label{kcl-3-opsd}
%		I_{\phi,i,j,t}=\mu_{\phi,i,j}I_{i,j,t}~\forall \phi,\forall i,\forall j\in \mathcal{X}_i,\forall t\\
		\label{kcl-4-opsd}
		\alpha_{\phi,i,j,t}\in \{0,1\}~\forall \phi,\forall i,\forall j,\forall t\\
		\label{kcl-5-opsd}
		\sum\nolimits_\phi{\alpha_{\phi,i,j,t}=1}~\forall i,\forall j,\forall t
		\end{eqnarray}
	\end{subequations}
	where
	$\alpha_{\phi,i,j,t}$ is a binary variable indicating whether the $j^\text{th}$ customer is connected to phase $\phi$ of node $i$ for period $t$ ($\alpha_{\phi,i,j,t}$ exists only when the corresponding customer is adjustable); 
	$\mu_{\phi,i,j}$ is a known parameter indicating the initial phase-position of $j^\text{th}$ customer connected to node $i$, with $\mu_{\phi,i,j}=1$ implying it is initially connected to phase $\phi$.
	
	%In the above expressions, $\alpha_{\phi,i,j,t}=1$ or $\mu_{\phi,i,j}=1$ means the phase position of the $j^\text{th}$ customer connected to node $i$ at period is $\phi$. Moreover, $\alpha_{\phi,i,j,t}$ is only introduced for each adjustable customer.
	
	\subsubsection{Constraints on service line for each customer}
	Ohm's law should be satisfied for each service line and the terminal voltages of PSD should be identical, leading to \eqref{kcl-6-opsd} and \eqref{kcl-7-opsd}, respectively.
	\begin{subequations}
		\label{kcl-opsd2}
		\begin{eqnarray}
		\label{kcl-6-opsd}
		U_{i,j,t}-V_{i,j,t}=Z_{i,j}I_{i,j,t}~\forall i,\forall j\in\mathcal{C}_i,\forall t\\
		\label{kcl-7-opsd}
		U_{i,j,t}=\sum\nolimits_{\phi}{\varepsilon_{\phi,i,j,t}V_{\phi,i,t}}~\forall i,\forall j,\forall t
%		\label{kcl-8-opsd}
%		U_{i,j,t}=\sum\nolimits_{\phi}{\mu_{\phi,i,j}V_{\phi,i,t}}~\forall i,\forall j\in\mathcal{X}_i,\forall t
		\end{eqnarray}
	\end{subequations}
	
	Obviously, \eqref{kcl-6-opsd} is a linear expression. However, \eqref{kcl-7-opsd} is linear for fixed customers and is non-convex for adjustable customers. 
	
	\subsubsection{Power balance equations}
	All customers are assumed to be with constant active and reactive powers, which yields
		\begin{eqnarray}
		\label{pb-1-opsd}
    	I_{i,j,t}=(P^n_{i,j,t}-jQ^n_{i,j,t})/V_{i,j,t}^*~\forall i,\forall j\in\mathcal{C}_i,\forall t
		\end{eqnarray}
	where $P^n_{i,j,t}$, $Q^n_{i,j,t}$ are the net active and reactive demands of $j^\text{th}$ customer at node $i$ for period $t$. 
	
	%It is noteworthy that $P^{g}_{\phi,i,j,t},Q^g_{\phi,i,j,t}$ are introduced for all fixed customers for clearer illustration. However, for customers without PV installed, $P^{g}_{\phi,i,j,t},Q^g_{\phi,i,j,t}$ should be removed.
	
	\subsubsection{Voltage constraints}
    The voltage magnitudes of all nodes and customers should always be within their limits, leading to
	\begin{subequations}
		\label{vlimit-opsd}
		\begin{eqnarray}
		\label{vlimit-1-opsd}
		V_{\phi,i}^\text{min}\le |V_{\phi,i,t}|\le V_{\phi,i}^\text{max}~\forall \phi,\forall i,\forall t\\
		\label{vlimit-2-opsd}
		V_{i,j}^\text{min}\le |V_{i,j,t}|\le V_{i,j}^\text{max}~\forall i,\forall j,\forall t
		\end{eqnarray}
	\end{subequations}
	where $V_{\phi,i}^\text{min}/V_{\phi,i}^\text{max}$ is the lower/upper voltage magnitude (VM) limit of $V_{\phi,i,t}$ and $V_{i,j}^\text{min}/V_{i,j}^\text{max}$ is the lower/upper VM limit of $V_{i,j,t}$.
	
	Particularly, the voltage of root node is assumed to be known as $V_{\phi,t}^0$, leading to
	\begin{eqnarray}
	\label{vroot-opsd}
	V_{\phi,x,t}=V_{\phi,t}^0=|V_{\phi,t}^0|\angle{\beta^0_{\phi,t}},\forall t
	\end{eqnarray}
	where $\beta^0_{\phi,t}$ is the phase angle of $V_{\phi,t}^0$.
	
	With known $V_{\phi,t}^0(\forall \phi)$, the phase angles of phase-to-phase voltages $V_{ab,t}^0,V_{bc,t}^0$ and $V_{ca,t}^0$ can be easily calculated. They are denoted as $\beta^0_{ab,t},\beta^0_{bc,t}$ and $\beta^0_{ca,t}$ for the convenience of later discussion. 

	\subsubsection{DT capacity limit}
	The currents of three phases running through the DT should not exceed their limits, leading to %capacity of Dt should not  of each line should not exceed its limit, leading to  
%	\begin{subequations}
%		\label{Ilimit-opsd}
		\begin{eqnarray}
		\label{Ilimit-1-opsd}
		|I_{\phi,xy,t}|\le I_{\phi,xy}^\text{max}~\forall \phi,\forall t
%		\label{Ilimit-1-PCD-opsd}
%		|I_{\phi,xy,t}-I_{\phi,t}|\le I_{\phi,xy}^\text{max}~\forall \phi,\forall t\\
%		\label{Ilimit-2-opsd}
%		|I_{i,j,t}|\le I_{i,j}^\text{max}~\forall i,\forall j\in\mathcal{C}_i,\forall t
		\end{eqnarray}
%	\end{subequations}
	where $I_{\phi,xy}^\text{max}$ is the upper limit for the magnitude of $I_{\phi,xy,t}$.% and $I_{i,j}^\text{max}$ is the upper CM limit of $I_{i,j,t}$.

	\subsubsection{Voltage unbalance constraints}
    Voltage unbalances, which can be defined by various methods \cite{RN51,NER,bess-ref-11}, should not exceed the specified level to improve the quality power supply. In this paper, the negative sequence voltage (NSV) and zero-sequence voltage (ZSV) for each node are constrained \cite{bess-ref-11}, yielding
	\begin{subequations}
		\label{negV_level-opsd}
		\small
		\begin{eqnarray}
		\label{negV-1-opsd}
		6V^-_{i,t}=2V_{a,i,t}-(1+j\sqrt{3})V_{b,i,t}-(1-j\sqrt{3})V_{c,i,t}~\forall i,\forall t\\
		\label{zsv-1-opsd}
		3V^0_{i,t}=V_{a,i,t}+V_{b,i,t}+V_{c,i,t}~\forall i,\forall t\\
		\label{negV-2-opsd}
		|V^-_{i,t}|\le \nu^-|V_{n}|,~|V^0_{i,t}|\le \nu^0|V_{n}|~\forall i,\forall t
		\end{eqnarray}
	\end{subequations}
	where $V^-_{i,t}$ and $V^0_{i,t}$ are the NSV and ZSV of node $i$ at period $t$, respectively; $\nu^-$ and $\nu^0$ are the coefficients for constraining their magnitudes; $|V_{n}|$ is the magnitude of the nominal voltage in the network, which equals to 1~p.u. in per unit system.% is the nominal voltage of LVDN.
	
	\subsubsection{Operational frequency of PSDs} 
    For a practical system, the phase positions of PSDs may be adjusted one or several times during the whole dispatching period to better mitigate the unbalances. Denoting the allowed times to adjust PSDs as $N_o$, $\mathcal{T}$ is divided into $N_o$ subsets evenly first, say $\mathcal{T}_1,\cdots,\mathcal{T}_k,\cdots,\mathcal{T}_{N_o}$. The operational frequency constraints for $\mathcal{T}_k$ can accordingly be formulated as
	\begin{eqnarray}
	\label{ofreq-opsd}
	\alpha_{\phi,i,j,t_1}=\alpha_{\phi,i,j,t_2}~\forall \phi, \forall i,\forall j,\forall t_1\in\mathcal{T}_k,\forall t_2\in\mathcal{T}_k
	\end{eqnarray}
	
%===========================================================================================
\subsection{Optimization problem in summary}
Based on the previously discussed objective function and operational constraints, the optimization problem is summarized as
\begin{eqnarray}
\small
\label{opsd-1}
\text{OPSD}:~\min\{\eqref{obj-1-opsd}|s.t.~\eqref{Ipcd-1}-\eqref{Ipcd-7},\eqref{ulevel-PCD-opsd}-\eqref{negV_level-opsd}~\forall t\in\mathcal{T}; \eqref{ofreq-opsd}~\forall k\}\nonumber
\end{eqnarray}

OPSD belongs to the challenging MINCP due to introduced integer variables, and non-convex constraints in \eqref{Ipcd-4}, \eqref{kcl-2-opsd}, \eqref{kcl-7-opsd}, \eqref{pb-1-opsd} and \eqref{vlimit-opsd}. Moreover, noting that the objective function is separable and only operational constraints belonging to the same dispatching period $\mathcal{T}_k$ are coupled via \eqref{ofreq-opsd}, the problem OPSD is equivalent to solving $N_o$ sub-problems, where the $k^{\text{th}}$ sub-problem can be formulated as 
\begin{eqnarray}
\small
\label{opsd-k}
\text{OPSD}_k:~\min\{\sum_{t\in\mathcal{T}_k}{z_{t}}|s.t.~\eqref{Ipcd-1}-\eqref{Ipcd-7},\eqref{ulevel-PCD-opsd}-\eqref{negV_level-opsd}~\forall t\in\mathcal{T}_k;\eqref{ofreq-opsd}\}\nonumber
\end{eqnarray}

%===========================================================================================
\section{Solution Techniques}\label{soltech-opsd}
\begin{figure}[htb!]
	\centering\includegraphics[scale=0.12]{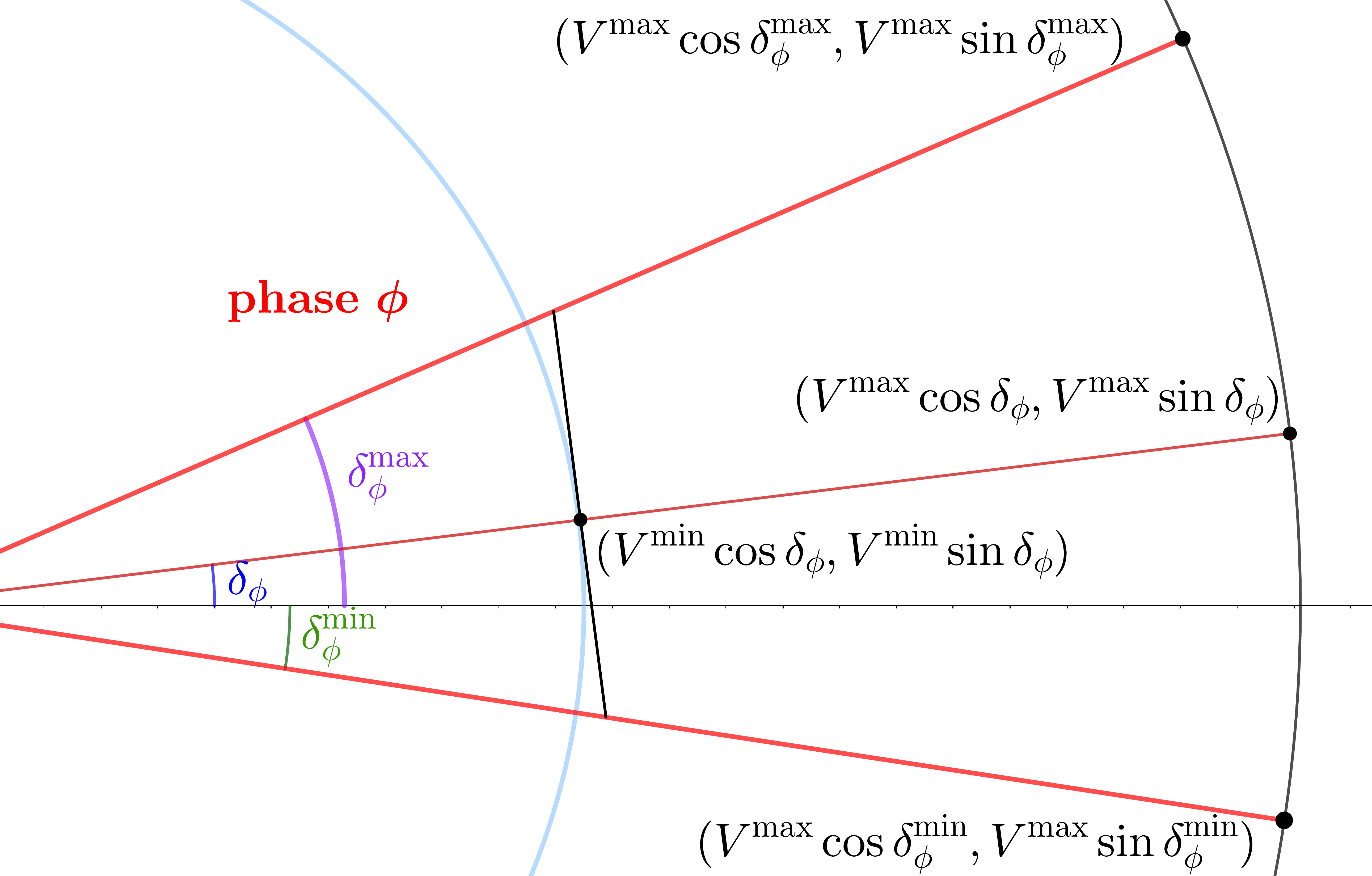}
	\caption{Illustration of linearizing voltage constraints}
	\label{fig-vtg}
\end{figure}

As OPSD belongs to MINCP that is difficult to solve, solution techniques to be discussed in this section aims at reformulating the problem to an efficient solvable problem. The techniques are based on the assumption that the difference of voltage angles in each phase is sufficiently small, which has been demonstrated in \cite{RN38,RN50,RN66,RN30,RN67}. As shown in Fig.\ref{fig-vtg}, VM limits at node $i$ for all phases are assumed to be within $[V^\text{min}_{i},V^\text{max}_{i}]$, and VA limits in phase $\phi$ for all nodes are within $[\delta^\text{min}_\phi,\delta^\text{max}_\phi]$. Further, VA for phase $\phi$ is assumed to be centered at $\delta_\phi$ and varies in $[\delta_\phi-\Delta\delta,\delta_\phi+\Delta\delta]$, which leads to $\delta^\text{min}_\phi=\delta_\phi-\Delta \delta,\delta^\text{max}_\phi=\delta_\phi+\Delta\delta$.  

The non-convexity of \eqref{Ipcd-4} is because $\beta_{\phi\psi,y,t}$ is a variable. However, with the assumption that the difference of voltage angles in each phase is sufficiently small, we further assume that $\beta_{\phi\psi,y,t}$ approximately equals to $\beta_{\phi\psi,x,t}=\beta^0_{\phi\psi,t}$. Then, \eqref{Ipcd-4} and \eqref{Ipcd-6} can be reformulated as
\begin{eqnarray}
\label{LIpcd-1}
I_{\phi\psi,t}=|I_{\phi\psi,t}|[(1-\kappa_{\phi\psi,t})\sin{\beta^0_{\phi\psi,t}}-\kappa_{\phi\psi,t}\sin{\beta^0_{\phi\psi,t}}\nonumber\\
+j\kappa_{\phi\psi,t}\cos{\beta^0_{\phi\psi,t}}-j(1-\kappa_{\phi\psi,t})\cos{\beta^0_{\phi\psi,t}}]~\forall\phi\psi\in\Phi,\forall t
\end{eqnarray}

Noting that the only non-convex part in \eqref{LIpcd-1} is $|I_{\phi\psi,t}|\kappa_{\phi\psi,t}$, which is the product of a continuous variable and a binary variable, it can be exactly reformulated as a mixed-integer linear (MIL) constraint noting the following equivalence. 
\begin{eqnarray}
\label{RLP-1}
\left.\begin{array}{r}
z=xy\\
x\in \{0,1\}\\
y\in[y^\text{min},y^\text{max}]
\end{array}
\right\}
\Leftrightarrow 
\left\{\begin{array}{r}
xy^\text{min}\le z \le xy^\text{max}\\
(x-1)y^\text{max}\le z-y\\
\le(x-1)y^\text{min}
\end{array}
\right.\
\end{eqnarray} 

For \eqref{kcl-2-opsd} and \eqref{kcl-7-opsd}, each of them for adjustable customers is also the product of a continuous and a binary variable, they can also be reformulated as MIL constraints based on \eqref{RLP-1}.

To reformulate power balance equation, \eqref{pb-1-opsd} is firstly linearized based on approximating $1/V^*$ as \eqref{fxy-opsd} \cite{RN38}.
\begin{eqnarray}
\label{fxy-opsd}
1/V^*_\phi\approx k^X_\phi X+k^Y_\phi Y+b^X_\phi+j(h^X_\phi X+h^Y_\phi Y+b^Y_\phi)
\end{eqnarray}
where $k^X_\phi,k^Y_\phi,b^X_\phi$ and $h^X_\phi,h^Y_\phi,b^Y_\phi$ are parameters to be fitted for phase $\phi$, which can be realized via least-square method (LSM) as discussed in \cite{RN38}. 

Moreover, the required known points to fit the parameters in \eqref{fxy-opsd} can be sampled in the feasible region as specified in Fig.\ref{fig-vtg}. With \eqref{fxy-opsd}, \eqref{pb-1-opsd} can be expressed as
\begin{subequations}
	\small
\label{pb-lin}
\begin{eqnarray}
\label{pb-1}
I_{i,j,t}=(P^n_{i,j,t}-jQ^n_{i,j,t})[k^X_{i,j,t}X_{i,j,t}+k^Y_{i,j,t}Y_{i,j,t}+b^X_{i,j,t}\nonumber\\
+j(h^X_{i,j,t}X_{i,j,t}+h^Y_{i,j,t}Y_{i,j,t}+b^Y_{i,j,t})]~\forall i,\forall j,\forall t\\
\label{pb-2}
k^X_{i,j,t}=\sum\nolimits_\phi{\varepsilon_{\phi,i,j,t}k^X_{\phi}},~h^X_{i,j,t}=\sum\nolimits_\phi{\varepsilon_{\phi,i,j,t}h^X_{\phi}}\\
\label{pb-3}
k^Y_{i,j,t}=\sum\nolimits_\phi{\varepsilon_{\phi,i,j,t}k^Y_{\phi}},~h^Y_{i,j,t}=\sum\nolimits_\phi{\varepsilon_{\phi,i,j,t}h^Y_{\phi}}\\
\label{pb-4}
b^X_{i,j,t}=\sum\nolimits_\phi{\varepsilon_{\phi,i,j,t}b^X_{\phi}},~b^Y_{i,j,t}=\sum\nolimits_\phi{\varepsilon_{\phi,i,j,t}b^Y_{\phi}}
\varepsilon
\end{eqnarray}
\end{subequations}

After replacing relevant terms in \eqref{pb-1} by \eqref{pb-2}-\eqref{pb-4}, non-convex parts will be introduced for adjustable customers, which are $z^C_{\phi,i,j,t}=\alpha_{\phi,i,j,t}X_{i,j,t}$ and $z^D_{\phi,i,j,t}=\alpha_{\phi,i,j,t}Y_{i,j,t}$ that can be further exactly reformulated as MIL constraints based on \eqref{RLP-1}.

The lower VM limits in \eqref{vlimit-opsd} are linearly approximated based on Fig.\ref{fig-vtg}, leading to 
\begin{subequations}
	\small
	\label{vm-lower}
	\begin{eqnarray}
	\label{vm-lower-1}
	X_{\phi,i,t}\cos{\delta_\phi}+Y_{\phi,i,t}\sin{\delta_\phi}\ge V^\text{min}_{\phi,i}~\forall \phi,\forall i,\forall t\\
	\label{vm-lower-2}
	X_{i,j,t}\sum\nolimits_\phi{\varepsilon_{\phi,i,j,t}\cos{\delta_\phi}}+Y_{i,j,t}\sum\nolimits_\phi{\varepsilon_{\phi,i,j,t}\sin{\delta_\phi}}\ge V^\text{min}_{i,j}\nonumber\\
	\forall i,\forall j,\forall t
	\end{eqnarray}
\end{subequations}

Obviously, \eqref{vm-lower-1} and \eqref{vm-lower-2} for fixed customers are linear. However, for adjustable customers, non-convex parts will be introduced after expanding \eqref{vm-lower-2}, which are also $z^C_{\phi,i,j,t}$ and $z^D_{\phi,i,j,t}$. It is noteworthy that the approximations are different for node voltages in the main feeder and terminal voltages of adjustable customers. This is because the parameters used for linearization in \eqref{vm-lower-2} are interdependent on the phase position of the customer as shown in Fig.\ref{fig-vtg}.

Moreover, there are nonlinear but convex constraints in OPSD, i.e. \eqref{Ipcd-5}, \eqref{ulevel-PCD-opsd-1}, upper limits in \eqref{vlimit-opsd}, \eqref{Ilimit-1-opsd}, \eqref{negV-2-opsd} and \eqref{negV-2-opsd}. Actually, they can be generally expressed as $\sqrt{x^2+y^2}\le z$, which is a second-order cone constraints and can be easily dealt with by commercial solvers.% and further reformulated as the following second-order cone constraint.
%\begin{eqnarray}
%\label{circle-soc}
%||x,y||_2\le z
%\end{eqnarray}

Several remarks on the formulation and solution techniques are given below.
\begin{enumerate}
	\item Based on the discussed solution techniques, OPSD is reformulated as a MISOCP problem, which can be solved via commercial solvers such as Cplex \cite{cplex} and Gurobi \cite{gurobi}. All programs for OPSD are implemented in Matlab with Yalmip \cite{yalmip}, and solved by Cplex 12.9 on a desktop PC with Intel i7-6700 3.4 GHz CPU, 16 GB memory. %Moreover, the time limit for the solver is set as 60 minutes to guarantee a feasible solution.
	\item As the original MINCP problem is strongly non-convex, errors will be inevitable by reformulating the problem as a MISOCP. Errors are partly from linearizing the power balance equation via \eqref{fxy-opsd} with given parameters $V^\text{min}_{i},V^\text{max}_{i},\delta^\text{min}_\phi$ and $\delta^\text{max}_\phi$. In this paper, $V^\text{min}_{i}$ and $V^\text{max}_{i}$ are specified according to practical operational requirements, while $\delta^\text{min}_\phi$ and $\delta^\text{max}_\phi$ can be determined via power flow analysis based on historical data to improve the accuracy of the approximation.
	\item Errors in the reformulation also come from linearizing lower voltage magnitude limits in \eqref{vm-lower}. However, noting that the linearized feasible region is a subset of the original non-convex feasible region, the original constraints can be guaranteed feasible with the optimal solution reported by MISOCP. %Moreover, appropriately selecting $\delta^\text{min}_\phi$ and $\delta^\text{max}_\phi$ will also help diminish the introduced errors as shown in Fig.\ref{fig-vtg}.
\end{enumerate}

%===========================================================================================
\section{Case Study}
\subsection{Case setup}
\begin{figure}[htb]
	\centering\includegraphics[scale=0.53]{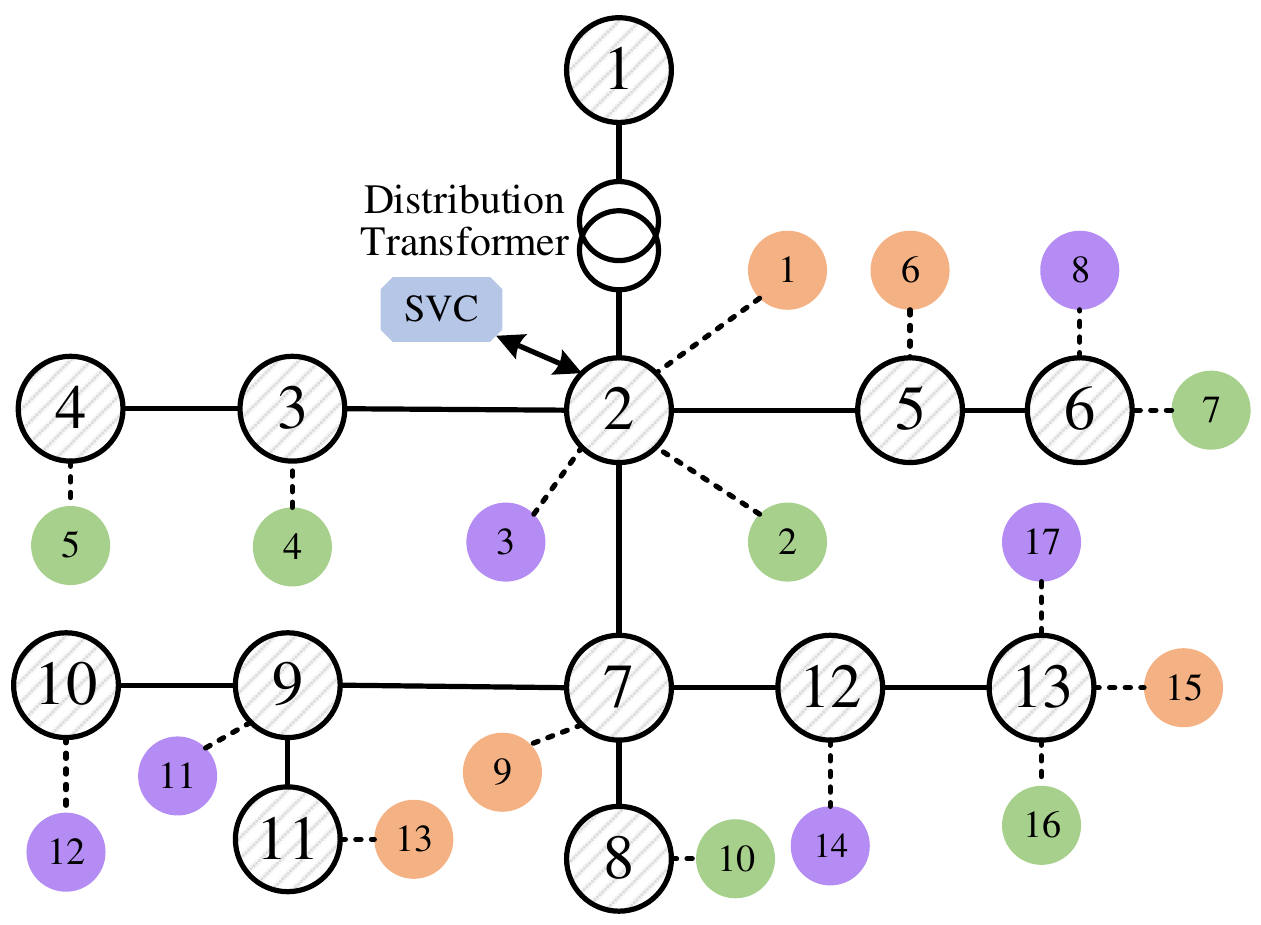}
	\caption{Topology of the modified IEEE-13 bus system (Orange: phase $a$; Green: phase $b$; Purple: phase $c$)}
	\label{fig-IEEE13-OPSD}
\end{figure}

Several cases based on the modified IEEE-13 bus system and a practical system in Australia will be studied. The topology of the modified IEEE-13 system is presented in Fig.\ref{fig-IEEE13-OPSD}. 

\begin{figure}[htb!]
	\centering\includegraphics[scale=0.18]{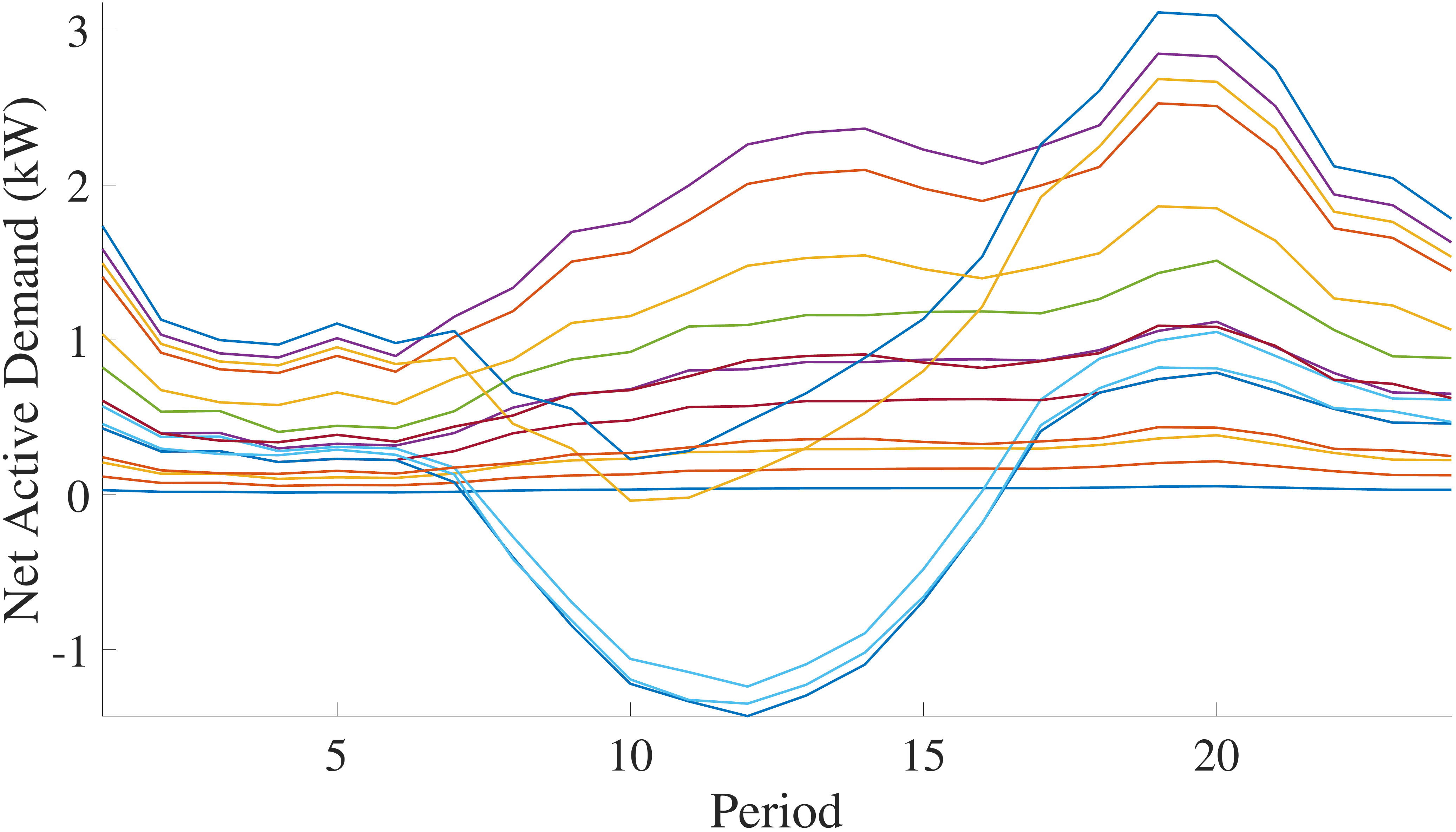}
	\caption{Net active demand of all residential customers for modified IEEE-13 system}
	\label{fig-dayload-IEEE13}
\end{figure}

The modified IEEE-13 system is with 17 single-phase powered customers, where 5 customers are with PV panels and the capacity for each of them is 2.5 kW. There are also 5 PSDs installed at customer 3,5,9,12 and 17. The net demands of all customers throughout a sunny day is presented in Fig.\ref{fig-dayload-IEEE13}, where PVs are assumed to be operating at their capacity values around the midday. 

The practical system in Australia is with 26 nodes in the main feeder and 69 customers that are powered with single-phase or three-phase. Among the 69 customers, 16 of them are with PV panels, which may further increase the unbalance levels when the solar irradiation levels are high around the midday. %The net demands of this system is presented in Fig.\ref{fig-dayload-Ausnet}. 

For all cases, $\nu^-$ is set as 2\% according to the operational requirements in Australian distribution networks \cite{AusCodes} and $\nu^0$ is set as 4.5\% in this paper. Other data or topology information of the two systems can be found in \cite{OPSD_data}.
%\begin{figure}[htb!]
%	\centering\includegraphics[scale=0.185]{fig-netP-Ausnet.pdf}
%	\caption{Net active demand of all residential customers for the system in Australia}
%	\label{fig-dayload-Ausnet}
%\end{figure}

The outline of the studied cases is given below.
\begin{enumerate}
	\item Case I: This case is based on the modified IEEE-13 bus system with all parameters specified in \cite{OPSD_data}. Moreover, $N_o$ is set as 1.    
	\item Case II: Same parameters as Case I are used for this case except that capacity for each PV panel is set as 3.0 kW.
	\item Case III: Parameters of this case are the same as Case I. However, the simulation results will be analyzed based on varying values of $N_o$.       
	\item Case IV: Parameters of this case are the same as Case I and the simulation results will be analyzed based on varying capacities of SVC.   
	\item Case V: This case is based on a practical system in Australia to test the practicality of the proposed method. 
\end{enumerate}

\subsection{Case I \& II}
\begin{figure}[htb!]
	\centering\includegraphics[scale=0.205]{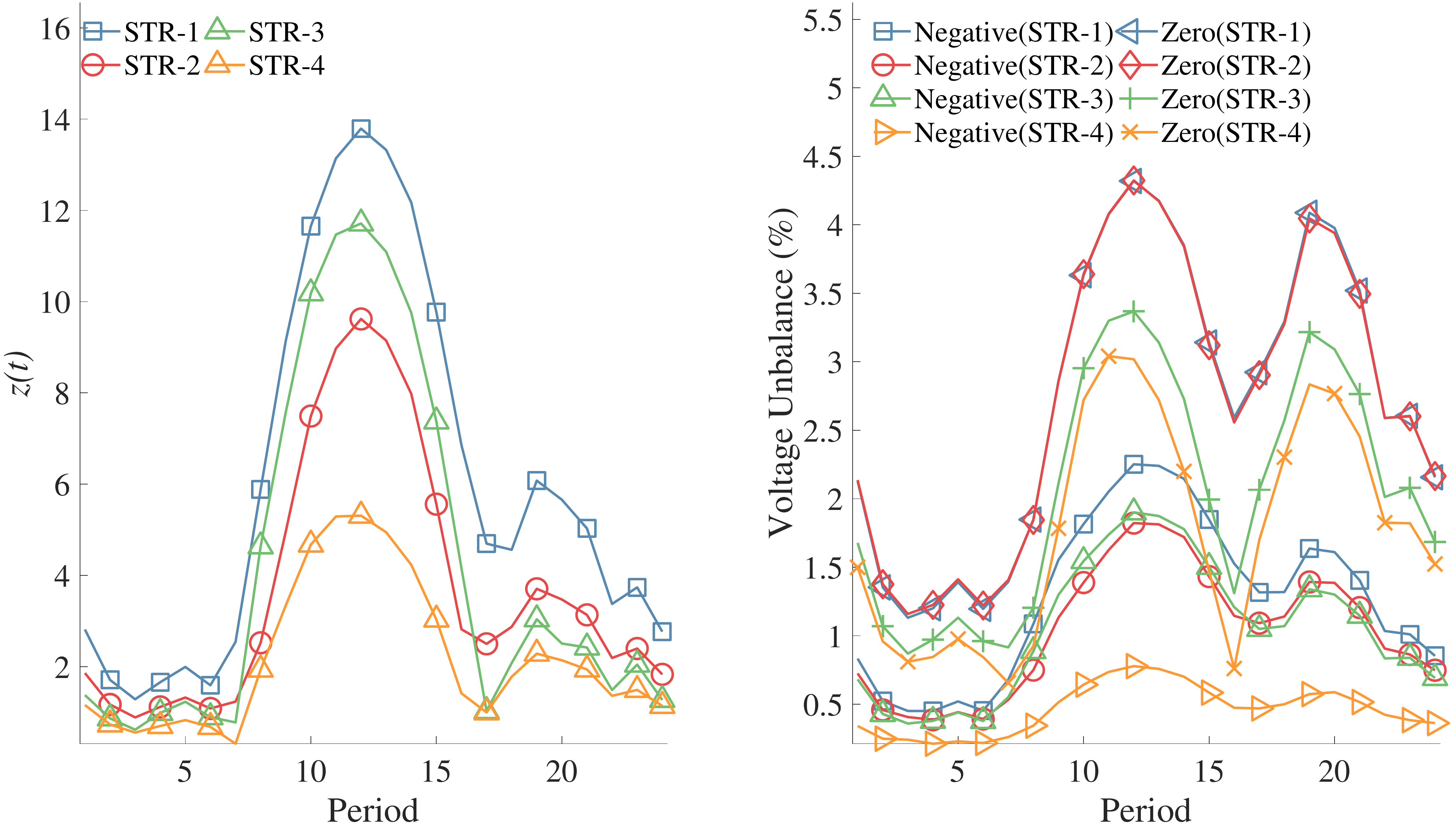}
	\caption{Current unbalance (left) and maximum voltage unbalance of all nodes (right) (STR-1: neither SVC nor PSD is used; STR-2: only SVC is used; STR-3: only PSD is used; STR-4: SVC and PSDs are used together)}
	\label{fig-OPSD-IEEE13-IandU-unbalance-SWTimes-01}
\end{figure}

In Case I, OPSD is successfully solved after 2.73s, 30.43s and 179.60s for STR-2, STR-3 and STR-4\footnote{The explanations for STR-1, STR-2, STR-3 and STR-4 are presented after the tile of Fig.\ref{fig-OPSD-IEEE13-IandU-unbalance-SWTimes-01}.}, respectively. When only PSD is used, the optimal phase positions for all PSDs are $b,b,c,a,a$, while the optimal results are $c,c,a,a,b$ when SVC and PSDs are used together. The unbalance level of currents running through the DT and the maximum voltage unbalance of all nodes in each period are presented in Fig.\ref{fig-OPSD-IEEE13-IandU-unbalance-SWTimes-01}.

When neither SVC nor PSD is employed in the network, strong current unbalances can be observed in the midday, which also results in the negative sequence voltage exceeding operational requirements. After introducing SVC or/and PSD, the negative sequence and zero-sequence voltage levels throughout the day all fall into required limits, which demonstrates that voltage unbalances can be well addressed by the two equipment. On current unbalances, employing either SVC or PSD can effectively reduce the overall unbalance level in the DT. However, the reduction will be much more significant if they are operated together as shown in the figure.  

\begin{figure}[htb!]
	\centering\includegraphics[scale=0.25]{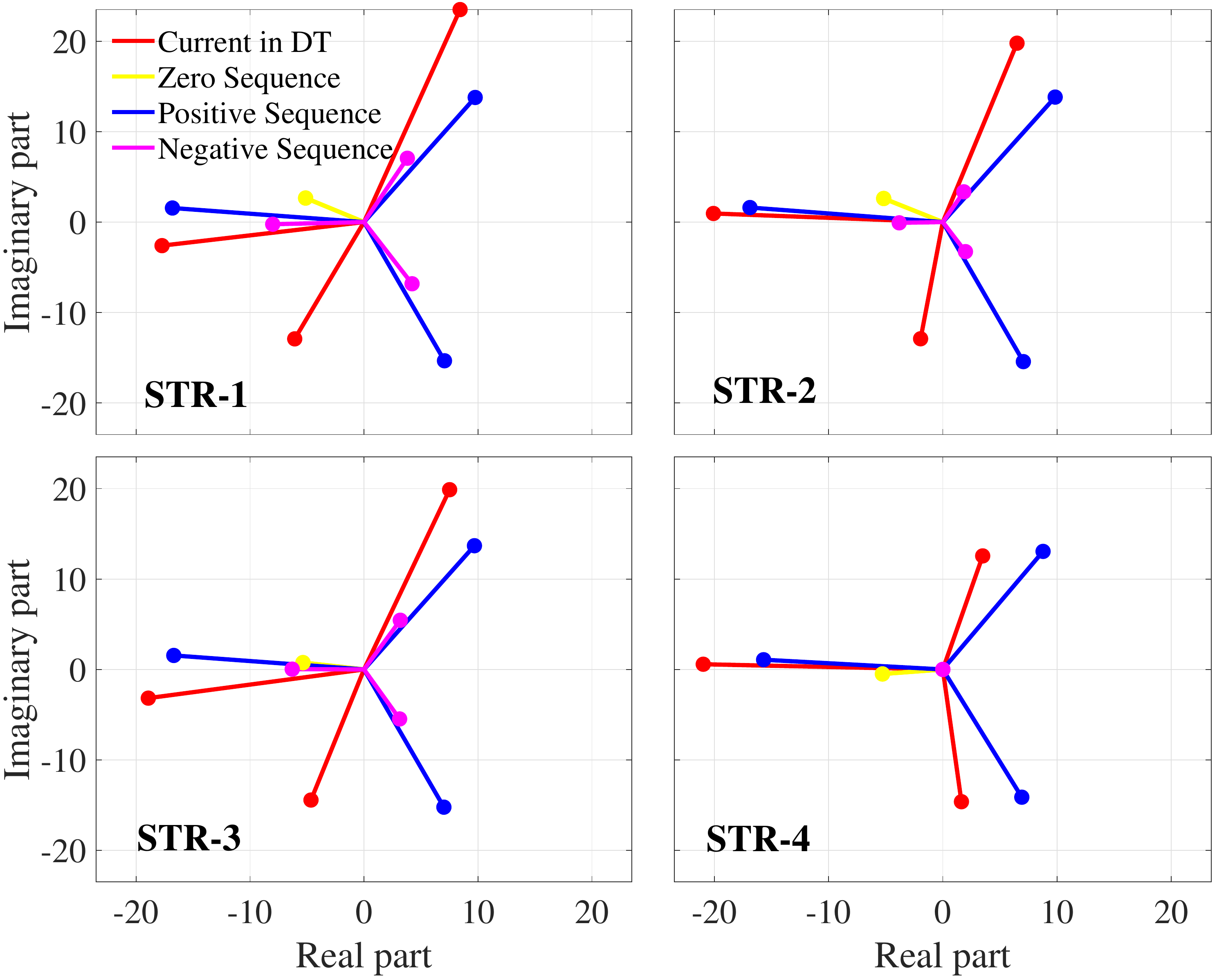}
	\caption{Decomposition of currents in the DT at period 12 for various operational strategies}
	\label{fig-decomp-case001}
\end{figure}

To show the current unbalance reduction, the currents of three phases running through the DT, and also their decomposition into zero, positive and negative sequences, at period 12 are presented in Fig.\ref{fig-decomp-case001}.

The simulation results clearly show that the strong current unbalances in the DT can be significantly mitigated if SVC and PSDs are used together. Moreover, using PSD and SVC together also outperforms the strategy when only one equipment is used, which demonstrates that more benefit can be brought to the network via using the two equipment together.  %As the values of $z_{14}$ for STR-2 and STR-3 are close to each other, it is interesting to note that the optimal solutions could be quite different, where STR-3 leads to more balanced currents in the DT. However, using SVC and PSDs together is highly recommended, which could bring the most benefit in addressing the unbalance issue.

\begin{figure}[htb!]
	\centering\includegraphics[scale=0.205]{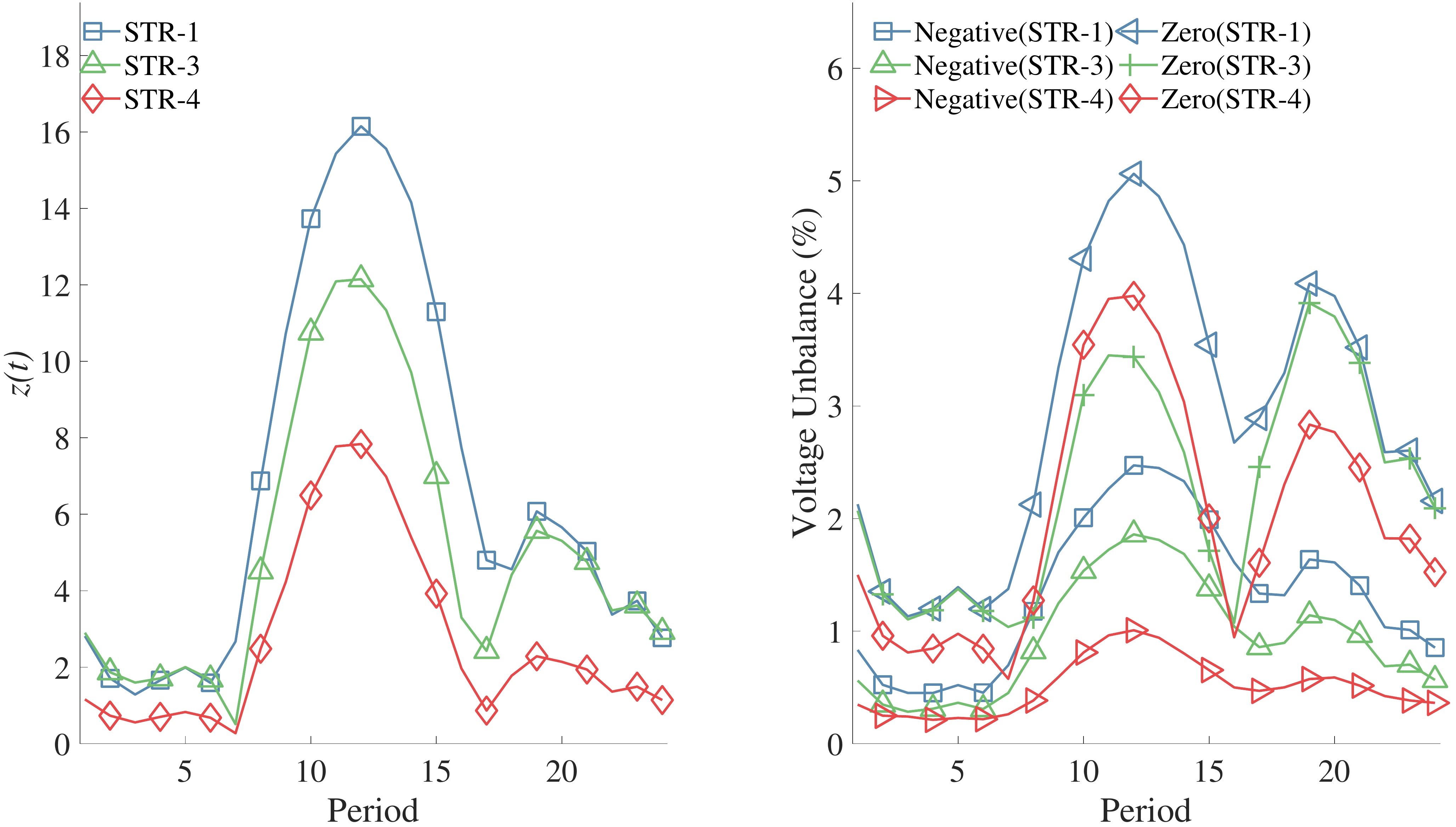}
	\caption{Current unbalance (left) and maximum voltage unbalance of all nodes (right) when each PV generation is 3.5 kW}
	\label{fig-OPSD-IEEE13-IandU-unbalance-SWTimes-01-PV3500}
\end{figure}

With higher PV penetration in Case II, infeasibility is reported for STR-2, which implies some operational requirements cannot be guaranteed by merely using SVC. When only PSD is used or SVC and PSDs are used together in this case, the optimal solutions are reported after 46.59s and 126.32s respectively, and the optimized current and voltage unbalances (STR-3 and STR-4) compared with the initial state (STR-1) are presented in Fig.\ref{fig-OPSD-IEEE13-IandU-unbalance-SWTimes-01-PV3500}.

Strong unbalances are also observed in the midday for STR-1, which leads to both negative sequence and zero-sequence voltages exceeding the operational requirements as shown in Fig.\ref{fig-OPSD-IEEE13-IandU-unbalance-SWTimes-01-PV3500}. However, by merely using PSDs or operating PSDs and SVC together, current unbalances in the midday can be significantly reduced while keeping the voltage unbalance levels under the upper limits. The simulation results also imply that the network's capability to balance the currents of three phases in the DT may be limited if only SVC is used and operating the two equipment together would bring more benefits.%It is noteworthy that although the overall unbalance level can be mitigated via STR-3, current unbalance levels may increase during peak hours in the night as shown in Fig.\ref{fig-OPSD-IEEE13-IandU-unbalance-SWTimes-01-PV3500}. However, this issue can be well addressed if both SVC and PSDs are operated together.%, which again demonstrates the potential benefit of operating them together. 

\subsection{Case III}
In practical operation, PSDs and SVC may be adjusted several times per day to achieve better operational performance. Specifically, we assume the control strategies will be updated every $T_o$ hours, where $T_o\in\{24,12,8,6,5,4\}$, thus leading to $N_o\in\{1,2,3,4,5,6\}$, respectively. As the system can be better balanced when SVC and PSDs are operated together, only STR-4 is studied in this section. The period subsets and computed current unbalances (objective values) under various $N_o$ are presented in Table \ref{tab-swtimes}.
\begin{table}[htb!]
	\footnotesize
	\caption{Period sets and objective values under various $N_o$}
	\label{tab-swtimes}
	\centering
	\setlength{\tabcolsep}{3pt}
	\renewcommand\arraystretch{0.9}
	\begin{tabular}{c|c|c}
		\hline\hline
		$N_o$& Period Subsets&\makecell{Objective Values (A)}\\
		\hline
		1&$\{1,\cdots,24\}$&52.2869\\	\hline
		2&$\{1,\cdots,12\},\{13,\cdots,24\}$&49.9326\\	\hline
		3&$\{1,\cdots,8\},\{9,\cdots,16\},\{17,\cdots,24\}$&16.0748\\	\hline
		4&\makecell{$\{1,\cdots,6\},\{7,\cdots,12\},\{13,\cdots,18\}$\\$\{19,\cdots,24\}$}&22.8034\\	\hline
		5&\makecell{$\{1,\cdots,5\},\{6,\cdots,10\},\{11,\cdots,15\}$\\$\{16,\cdots,20\},\{21,\cdots,24\}$}&17.5391\\	\hline
		6&\makecell{$\{1,\cdots,4\},\{5,\cdots,8\},\{9,\cdots,12\}$\\$\{13,\cdots,16\},\{17,\cdots,20\},\{21,\cdots,24\}$}&14.7673\\
		\hline\hline
	\end{tabular}
\end{table}

According to Table \ref{tab-swtimes} and the formulation of OPSD, it can be concluded that $F_{4}\le F_{2}\le F_{1}$, $F_6\le F_3\le F_1$ and $F_5\le F_1$, which can be well demonstrated by the results in Table \ref{tab-swtimes}. This is because the former one in each inequality has a larger feasible region than the latter one. However, it is interesting to note that for the studied case we have $F_4>F_3$ and $F_5>F_3$, implying that increasing allowed switching times does not always guarantee a reduction of unbalance levels. In other words, the timing to adjust PSDs is also important and needs to be carefully selected in practical application. 

\subsection{Case IV}
\begin{figure}[htb!]
	\centering\includegraphics[scale=0.205]{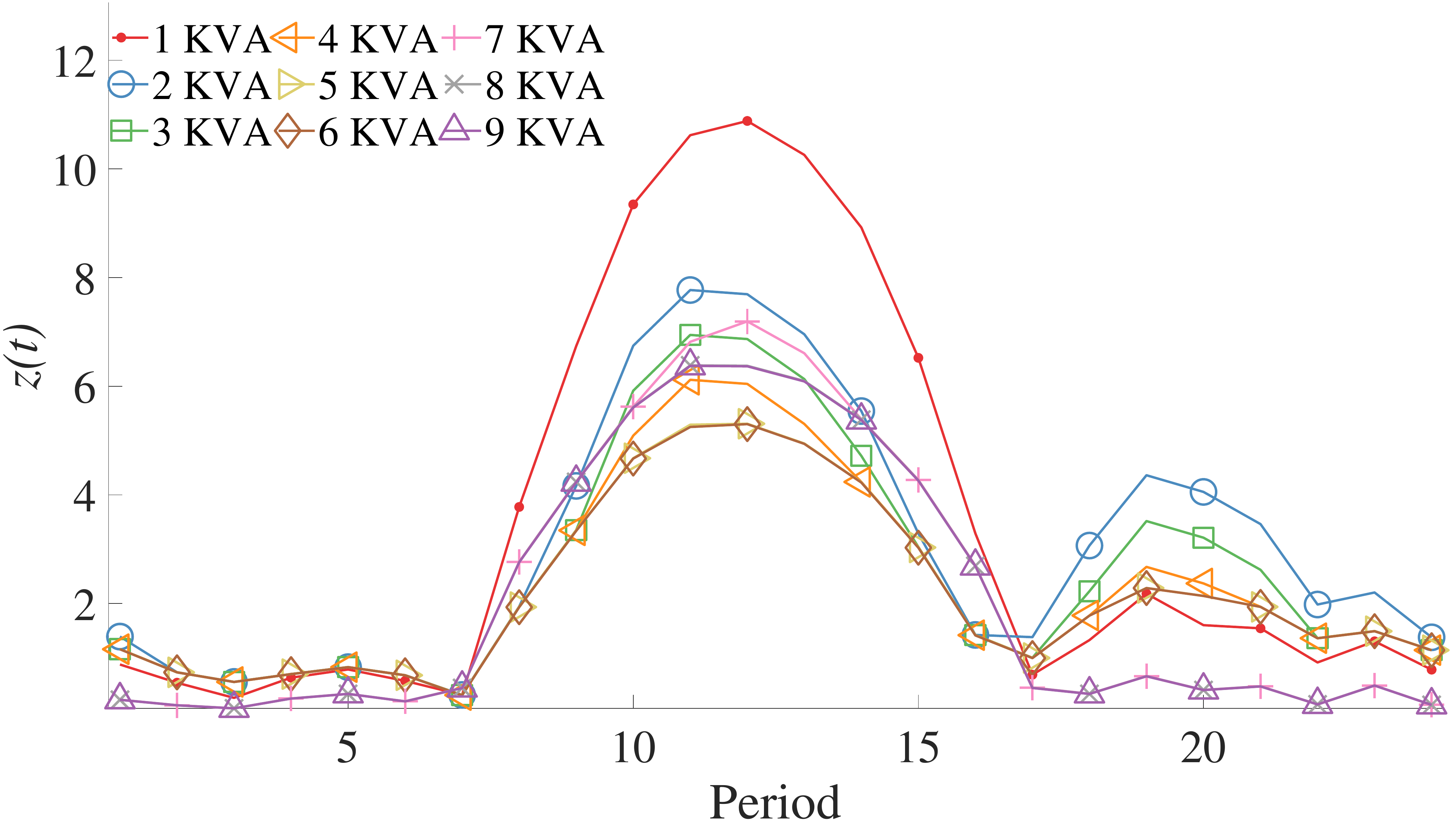}
	\caption{Current unbalances under various capacities of SVC (STR-4)}
	\label{fig-OPSD-IEEE13-diffSVG}
\end{figure}

The impact of SVC capacity on mitigating unbalances will be investigated in this case. Specifically, the capacity of SVC will be made to vary across a range of values $[1~\text{kVA},2~\text{kVA},\cdots,8~\text{kVA},9~\text{kVA}]$ and the OPSD will be solved successively with $N_o=1$ for STR-4. Simulation results are presented in Fig.\ref{fig-OPSD-IEEE13-diffSVG}.

Simulation results show that the objective value decreases with larger size of SVC quipped. $z_t$ for all periods are presented in Fig.\ref{fig-OPSD-IEEE13-diffSVG}, which shows that the marginal benefit in mitigating current balance is extremely small when SVC capacity exceeds 8 kVA. Therefore, With enough historic data and the targeted unbalance reduction to be achieved, the proposed model can be used to determine suitable capacity of SVC for a practical network. 

\subsection{Case V}
\begin{figure}[htb!]
	\centering\includegraphics[scale=0.205]{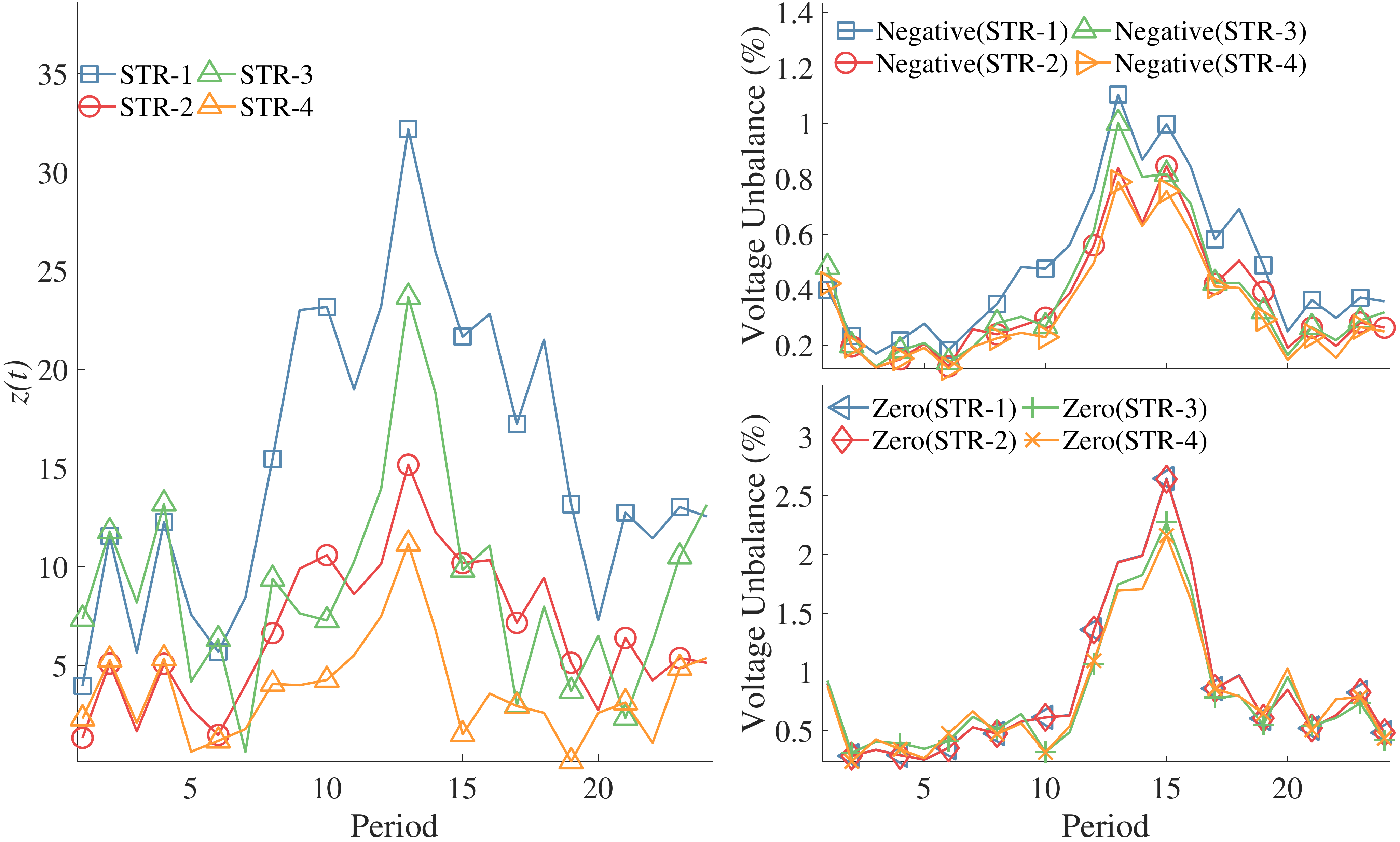}
	\caption{Current unbalance (left) and maximum voltage unbalance of all nodes (right)}
	\label{fig-OPSD-Ausnet-IandU-unbalance-SWTimes-01}
\end{figure}

\begin{figure}[htb!]
	\centering\includegraphics[scale=0.25]{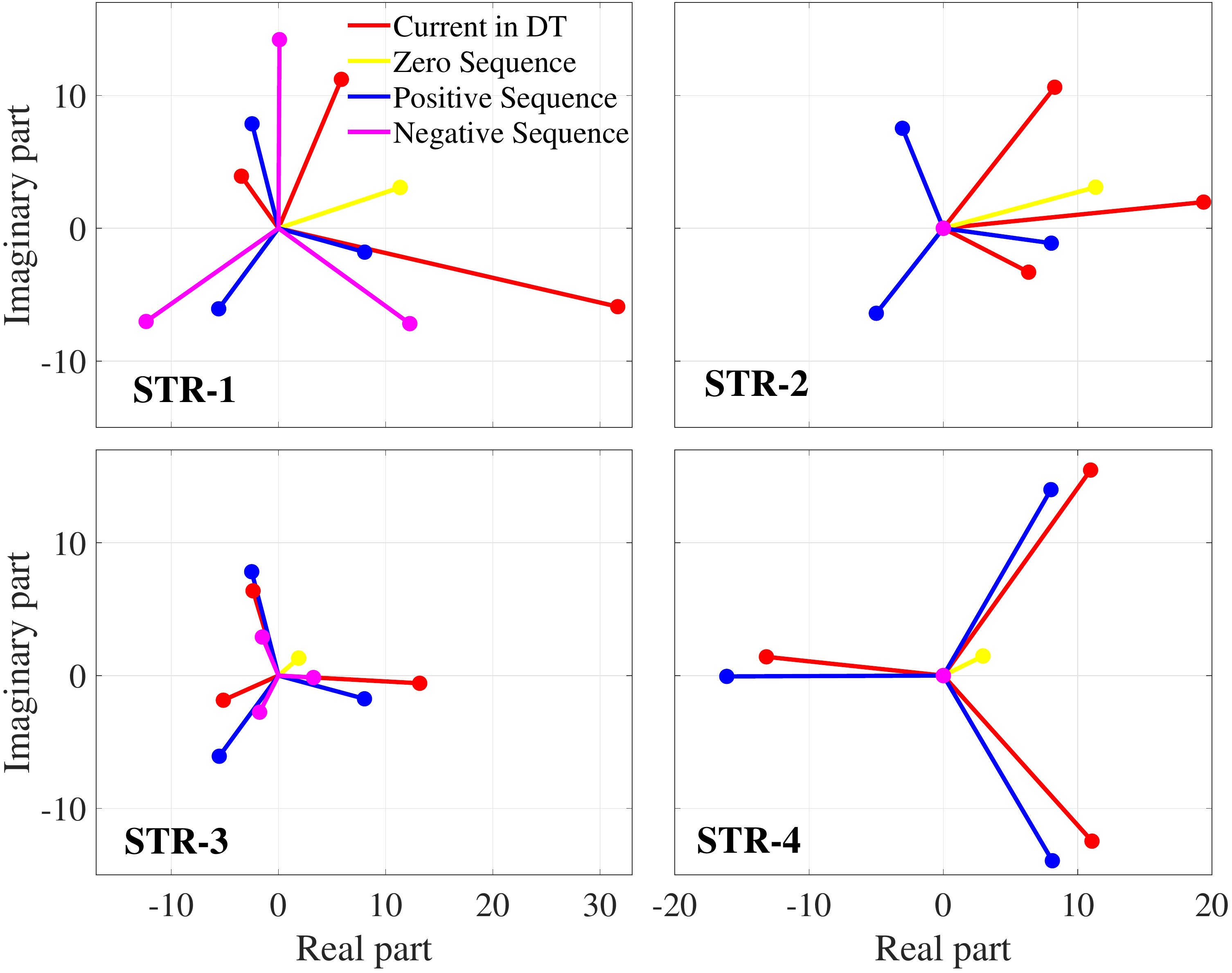}
	\caption{Decomposition of currents in the DT at period 14 for various operational strategies}
	\label{fig-decomp-case005}
\end{figure}

The current and voltage unbalances after optimization for Case V are presented in Fig.\ref{fig-OPSD-Ausnet-IandU-unbalance-SWTimes-01} and Fig.\ref{fig-decomp-case005} with $N_o=1$. Similar to the simulation results of Case I, current unbalance levels can be significantly reduced when SVC and PSDs are operated together in the network. It is noteworthy that although both negative sequence and zero-sequence voltage levels are within their upper limits, they are overall decreased in a more balanced network (STR-4).

Regarding the computational efficiency, the solver takes 8.84s, 1172.2s and 3718.9s\footnote{A sub-optimal solution that has a close objective value with the optimal solution was detected just after 403.09s.} to report the optimal solutions for STR-2, STR-3 and STR-4, respectively, which well demonstrates the practicality of the presented method to mitigate unbalances in LVDN via SVC and PSDs.

%===========================================================================================
\section{Conclusions}
Mitigating unbalances in LVDN via SVC and PSDs is formulated as a MINCP and solved by commercial solvers after reformulating it as a MISOCP. Based on the simulation results, major conclusions are:
\begin{enumerate}
	\item Employing SVC or/and PSD can effectively mitigate unbalances in LVDN. However, better performance can be achieved if they are operated together. Moreover, merely using one equipment sometimes cannot provide a feasible solution, which demonstrates the necessity of coordinately operating the two equipment in some cases.%using them together can help better address the unbalance issue and provide an secure operational strategy.
	\item Increasing allowed times to adjust PSDs during the whole dispatching period can help better address unbalance issues in the network. However, the improvement cannot be guaranteed because the timing to adjust the PSDs is also critical.
	\item The proposed optimization model can be used to determine the capacity of SVC in a practical network or assess whether the existing capacity of SVC is sufficient to achieve the target of unbalance reduction. 
\end{enumerate} 

Developing more efficient algorithms and incorporating other controllable resources, e.g. active/reactive powers from PVs and residential or network-based battery energy storage system (BESS), fall in our future research interests. 

%XXXXXX
%
%It is noteworthy that the voltage of root node is assumed to be known in the formulation. However, in practical system, the operation strategy of both PSDs and SVC will affect the power flows in the upstream MVDN, which will reversely influence voltage of root node in LVDN. How to coordinate the operation of MVDN an LVDN with these equipment falls in our future research interest.	
	
\bibliography{OPSDREF}
\end{document}